\magnification=\magstep1


\def\item{\vskip1.3pt\hang\textindent}


\tolerance=300 \pretolerance=200 \hfuzz=1pt \vfuzz=1pt

\hoffset 0cm            
\hsize=5.8 true in \vsize=9.5 true in

\def\rightheadline{\hfil\smc\lastname\hfil\tenbf\folio}
\def\leftheadline{\tenbf\folio\hfil\smc\lastname\hfil}
\headline={\ifodd\pageno\rightheadline\else\leftheadline\fi}
\newdimen\dimenone
\def\checkleftspace#1#2#3#4#5{
 \dimenone=\pagetotal
 \advance\dimenone by -\pageshrink   
 \ifdim\dimenone>\pagegoal          
   \else\dimenone=\pagetotal
        \advance\dimenone by \pagestretch
        \ifdim\dimenone<\pagegoal
          \dimenone=\pagetotal
          \advance\dimenone by#1         
          \setbox0=\vbox{#2\parskip=0pt                
                       \hyphenpenalty=10000
                       \rightskip=0pt plus 5em
                       \noindent#3 \vskip#4}    
        \advance\dimenone by\ht0
        \advance\dimenone by 3\baselineskip
        \ifdim\dimenone>\pagegoal\vfill\eject\fi
          \else\eject\fi\fi}

\parindent=35pt
\mathsurround=1pt
\parskip=1pt plus .25pt minus .25pt
\normallineskiplimit=.99pt

\mathchardef\emptyset="001F 

\def\Int{\mathop{\rm int}\nolimits}
%

\def\L{\mathop{\bf L}\nolimits}

\def\1{{\bf1}}\def\0{{\bf0}}

\def\({\bigl(}  \def\){\bigr)}
\def\<{\mathopen{\langle}}\def\>{\mathclose{\rangle}}

\def\Z{{\mathchoice{{\hbox{$\rm Z\hskip 0.26em\llap{\rm Z}$}}}%
{{\hbox{$\rm Z\hskip 0.26em\llap{\rm Z}$}}}%
{{\hbox{$\scriptstyle\rm Z\hskip 0.31em\llap{$\scriptstyle\rm Z$}$}}}{{%
\hbox{$\scriptscriptstyle\rm
Z$\hskip0.18em\llap{$\scriptscriptstyle\rm Z$}}}}}}

\def\N{{\mathchoice{\hbox{$\rm I\hskip-0.14em N$}}%
{\hbox{$\rm I\hskip-0.14em N$}}%
{\hbox{$\scriptstyle\rm I\hskip-0.14em N$}}%
{\hbox{$\scriptscriptstyle\rm I\hskip-0.10em N$}}}}

\def\R{{\mathchoice{\hbox{$\rm I\hskip-0.14em R$}}%
{\hbox{$\rm I\hskip-0.14em R$}}%
{\hbox{$\scriptstyle\rm I\hskip-0.14em R$}}%
{\hbox{$\scriptscriptstyle\rm I\hskip-0.10em R$}}}}

\def\K{{\mathchoice{\hbox{$\rm I\hskip-0.15em K$}}%
{\hbox{$\rm I\hskip-0.15em K$}}%
{\hbox{$\scriptstyle\rm I\hskip-0.15em K$}}%
{\hbox{$\scriptscriptstyle\rm I\hskip-0.11em K$}}}}

\def\qed{\hfill {\hbox{[\hskip-0.05em ]}}}

\def\.{{\cdot}}
\def\|{\Vert}
\def\ssk{\smallskip}
\def\msk{\medskip}
\def\bsk{\bigskip}
\def\giantskip{\vskip2\bigskipamount}

\def\giantbreak{\par \ifdim\lastskip<2\bigskipamount \removelastskip
         \penalty-400 \giantskip\fi}

\def\nin{\noindent}
\def\cen{\centerline}
\def\pagebreak{\vskip 0pt plus 0.0001fil\break}
\def\linebreak{\break}

\def\epsilon{\varepsilon}

\font\ninerm=cmr9 \font\eightrm=cmr8 \font\sixrm=cmr6
 \font\eightbf=cmbx8 \font\sixbf=cmbx6
 \font\eighti=cmmi8 \font\sixi=cmmi6
\font\ninesy=cmsy9 \font\eightsy=cmsy8 \font\sixsy=cmsy6
 \font\eightit=cmti8 
 \font\eightsl=cmsl8 
\font\eighttt=cmtt8
\font\bfone=cmbx10 scaled\magstep1 
\font\smc=cmcsc10 
 
scaled\magstep1 \font\small=cmcsc8

\def\no #1. {\bigbreak\vskip-\parskip\noindent\bf #1. \quad\rm}

\def\Proposition #1. {\checkleftspace{0pt}{\bf}{Theorem}{0pt}{}
\bigbreak\vskip-\parskip\noindent{\bf Proposition #1.} \quad\it}

\def\Theorem #1. {\checkleftspace{0pt}{\bf}{Theorem}{0pt}{}
\bigbreak\vskip-\parskip\noindent{\bf  Theorem #1.} \quad\it}
\def\Corollary #1. {\checkleftspace{0pt}{\bf}{Theorem}{0pt}{}
\bigbreak\vskip-\parskip\nin{\bf Corollary #1.} \quad\it}
\def\Lemma #1. {\checkleftspace{0pt}{\bf}{Theorem}{0pt}{}
\bigbreak\vskip-\parskip\noindent{\bf  Lemma #1.}\quad\it}

\def\Definition #1. {\checkleftspace{0pt}{\bf}{Theorem}{0pt}{}
\rm\bigbreak\vskip-\parskip\noindent{\bf Definition #1.} \quad}

\def\Remark #1. {\checkleftspace{0pt}{\bf}{Theorem}{0pt}{}
\rm\bigbreak\vskip-\parskip\noindent{\bf Remark #1.}\quad}

\def\Exercise #1. {\checkleftspace{0pt}{\bf}{Theorem}{0pt}{}
\rm\bigbreak\vskip-\parskip\noindent{\bf Exercise #1.} \quad}

\def\Example #1. {\checkleftspace{0pt}{\bf}{Theorem}{0pt}{}
\rm\bigbreak\vskip-\parskip\noindent{\bf Example #1.}\quad}
\def\Examples #1. {\checkleftspace{0pt}{\bf}{Theorem}{0pt}
\rm\bigbreak\vskip-\parskip\noindent{\bf Examples #1.}\quad}

\newcount\problemnumb \problemnumb=0
\def\Problem{\global\advance\problemnumb by 1\bigbreak\vskip-\parskip\noindent
{\bf Problem \the\problemnumb.}\quad\rm }

\def\Proof#1.{\rm\par\ifdim\lastskip<\bigskipamount\removelastskip\fi\smallskip
            \noindent {\bf Proof.}\quad}

\nopagenumbers

\def\author{}
\def\lastname{}
\def\thanks#1{\footnote*{\eightrm#1}}
\def\title{}

\def\nonumbers{\def\leftheadline{\hfil} \def\rightheadline{\hfil}}

\def\lastname{}
\def\h{{\textstyle{1\over2}}}

\def\n{{\cal N}}
\def\ep{\epsilon}

\def\text{\textstyle}
\def\disp{\displaystyle}
\def\d{{\,\rm d}}

\def\L{{\bf L\"osung:\quad}}

\def\and{{\rm and }}

\def\n{\cen{{\it W.G. Nowak}}}

\expandafter\edef\csname amssym.def\endcsname{%
       \catcode`\noexpand\@=\the\catcode`\@\space}
\catcode`\@=11
\def\undefine#1{\let#1\undefined}
\def\newsymbol#1#2#3#4#5{\let\next@\relax
 \ifnum#2=\@ne\let\next@\msafam@\else
 \ifnum#2=\tw@\let\next@\msbfam@\fi\fi
 \mathchardef#1="#3\next@#4#5}
\def\mathhexbox@#1#2#3{\relax
 \ifmmode\mathpalette{}{\m@th\mathchar"#1#2#3}%
 \else\leavevmode\hbox{$\m@th\mathchar"#1#2#3$}\fi}
\def\hexnumber@#1{\ifcase#1 0\or 1\or 2\or 3\or 4\or 5\or 6\or 7\or 8\or
 9\or A\or B\or C\or D\or E\or F\fi}


\font\tenmsb=msbm10 \font\sevenmsb=msbm7 \font\fivemsb=msbm5
\newfam\msbfam
\textfont\msbfam=\tenmsb\scriptfont\msbfam=\sevenmsb
\scriptscriptfont\msbfam=\fivemsb \edef\msbfam@{\hexnumber@\msbfam}
\def\Bbb#1{{\fam\msbfam\relax#1}}

\font\teneufm=eufm10 \font\seveneufm=eufm7 \font\fiveeufm=eufm5
\newfam\eufmfam
\textfont\eufmfam=\teneufm \scriptfont\eufmfam=\seveneufm
\scriptscriptfont\eufmfam=\fiveeufm

\catcode`@=11 

\expandafter\edef\csname amssym.def\endcsname{%
       \catcode`\noexpand\@=\the\catcode`\@\space}
\font\eightmsb=msbm8 \font\sixmsb=msbm6 \font\fivemsb=msbm5
\font\eighteufm=eufm8 \font\sixeufm=eufm6 \font\fiveeufm=eufm5
\newskip\ttglue
\def\eightpoint{\def\rm{\fam0\eightrm}%
  \textfont0=\eightrm \scriptfont0=\sixrm \scriptscriptfont0=\fiverm
  \textfont1=\eighti \scriptfont1=\sixi \scriptscriptfont1=\fivei
  \textfont2=\eightsy \scriptfont2=\sixsy \scriptscriptfont2=\fivesy
  \textfont3=\tenex \scriptfont3=\tenex \scriptscriptfont3=\tenex
\textfont\eufmfam=\eighteufm \scriptfont\eufmfam=\sixeufm
\scriptscriptfont\eufmfam=\fiveeufm \textfont\msbfam=\eightmsb
\scriptfont\msbfam=\sixmsb \scriptscriptfont\msbfam=\fivemsb
  \def\it{\fam\itfam\eightit}%
  \textfont\itfam=\eightit
  \def\sl{\fam\slfam\eightsl}%
  \textfont\slfam=\eightsl
  \def\bf{\fam\bffam\eightbf}%
  \textfont\bffam=\eightbf \scriptfont\bffam=\sixbf
   \scriptscriptfont\bffam=\fivebf
  \def\tt{\fam\ttfam\eighttt}%
  \textfont\ttfam=\eighttt
  \tt \ttglue=.5em plus.25em minus.15em
  \normalbaselineskip=9pt
  \def\MF{{\manual opqr}\-{\manual stuq}}%
  \let\big=\eightbig
  \setbox\strutbox=\hbox{\vrule height7pt depth2pt width\z@}%
  \normalbaselines\rm}
\def\eightbig#1{{\hbox{$\textfont0=\ninerm\textfont2=\ninesy
  \left#1\vbox to6.5pt{}\right.\n@space$}}}


\csname amssym.def\endcsname


\def\al{\alpha}
\def\be{\beta}

\def\({\left(}
\def\){\right)}

\def\eq{\eqalign}

\def\O#1{O\(#1\)}
\def\abs#1{\left| #1 \right|}

\def\klein{\eightpoint \def\smc{\small} \baselineskip=9pt}

\def\fn#1#2{{\parindent=0.7true cm
\footnote{$^{(#1)}$}{{\klein  #2}}}}

\font\boldmas=msbm10                  
\def\Bbb#1{\hbox{\boldmas #1}}        
\def\Z{{\Bbb Z}}                        
\def\N{{\Bbb N}}                        

\def\R{{\Bbb R}}


\font\eightrm=cmr8 \long\def\fussnote#1#2{{\baselineskip=9pt
\setbox\strutbox=\hbox{\vrule height 7pt depth 2pt width 0pt}%
\eightrm \footnote{#1}{#2}}}
\font\boldmasi=msbm10 scaled 700      
\def\Bbbi#1{\hbox{\boldmasi #1}}      
\font\boldmas=msbm10                  
\def\Bbb#1{\hbox{\boldmas #1}}        
\def\Zi{{\Bbbi Z}}                      
\def\Pi{{\Bbbi P}}                      



\def\dint #1 {
\quad  \setbox0=\hbox{$\disp\int\!\!\!\int$}
  \setbox1=\hbox{$\!\!\!_{#1}$}
  \vtop{\hsize=\wd1\centerline{\copy0}\copy1} \quad}

\def\drint #1 {
\qquad  \setbox0=\hbox{$\disp\int\!\!\!\int\!\!\!\int$}
  \setbox1=\hbox{$\!\!\!_{#1}$}
  \vtop{\hsize=\wd1\centerline{\copy0}\copy1}\qquad}

\def\frac#1#2{{#1\over #2}}

\def\date{\the\day.~\the\month.~\the\year}

\def\klein{\eightpoint \def\smc{\small} }

\def\frac#1#2{{#1\over#2}}
\def\Int{\int\limits}

\def\vol{{\rm vol}}

\nonumbers

\hsize=16true cm     \vsize=24true cm

\parindent=0cm

\def\b#1{{\bf #1}}
\def\nm#1{\left|#1\right|_*}
\def\vol{{4\pi\over3}}
\def\z{{[2]}}
\def\wz{\sqrt{2}\,}
\def\n{{N<n\le U}}
\def\m{{M-H<m\le W}}
\def\mk{M<m+k\le W}
\def\wt{\sqrt{t}}
\def\E{{\cal E}}
\def\K{{\cal K}}
\def\B{{\cal B}}
\def\L{{\cal L}}
\def\I{^{({\rm I})}}
\def\II{^{({\rm II})}}
\def\III{^{({\rm III})}}

\vbox{\vskip1true cm}

\cen{\bfone Eine effektive Absch\"{a}tzung f\"{u}r die } \msk
\cen{\bfone Gitter-Diskrepanz von Rotationsellipsoiden}\bsk

\cen{{\bf Ekkehard Kr\"{a}tzel und Werner Georg Nowak}\fn{*}{Die
Autoren danken dem \"{O}sterreichischen Fonds zur F\"{o}rderung der
wissenschaftlichen Forschung (FWF) f\"{u}r finanzielle
Unterst\"{u}tzung unter der Projekt-Nr.~P18079-N12.} {\bf(Wien)}}
\bsk\bsk

\vbox{\vskip 1true cm}

\footnote{}{\ssk \klein{\it Mathematics Subject Classification }
(2000): 11P21, 11K38, 52C07.\par }

{\klein{\bf Abstract. \ An effective estimate for the lattice point
discrepancy of ellipsoids of rotation. } For the lattice point
discrepancy $P_\E(x)$ (i.e., the number of integer points minus the
volume) of the ellipsoid $ (u_1^2+u_2^2)/a+a^2\,u_3^2\le x$ \ ($a,
x>0$), this paper provides an estimate of the shape \ssk
\cen{$|P_\E(x)|\le1237\,a^{1/8}\,x^{11/16}\,(\log(100x)+|\log
a|)^{3/8}\ +$ terms of smaller order in $x$.} }

\vbox{\vskip 1.2true cm}

{\bf 1.~Einleitung. } Die Theorie der Gitterpunkte in "gro{\ss}en"
Bereichen (im Sinne von E.~Landau) besitzt eine lange und sehr
erfolgreiche Geschichte, die z.B.~von E.~Kr\"{a}tzel in den
Monographien [8], [9] umfassend dargestellt wurde. Man vgl.~dazu
auch einen neueren \"{U}bersichtsartikel der Autoren mit A.~Ivi\'c
und M.~K\"{u}hleitner [7]. Kernproblem dieser Theorie ist die
Absch\"{a}tzung der Gitter-Diskrepanz (Gitterrest)
$$ P_\K(x)  = \#\(\Z^s\cap \sqrt{x}\,\K\) - \hbox{vol}(\K)x^{s/2} \eqno(1.1)
$$ f\"{u}r die um einen Faktor $\sqrt
x$ linear vergr\"{o}{\ss}erte Kopie $\sqrt{x}\,\K$ eines allgemeinen
oder spezielleren Bereiches $\K$ der Ebene und der R\"{a}ume $\R^s$,
$s\ge3$. Die Resultate wurden durchwegs in der Landau'schen
$O$-Notation formuliert. \ssk Der ebenso naheliegende wie
berechtigte Gedanke, die darin jeweils involvierten Konstanten unter
Kontrolle zu bekommen, taucht nur einerseits sehr fr\"{u}h bei
J.G.~Van der Corput [14] auf, andererseits erst in neuester Zeit: So
bewiesen V.~Bentkus und F.~G\"otze [1] bzw. F.~G\"{o}tze [4] f\"{u}r
ein beliebiges $\b o$-symmetrisches Ellipsoid $\E$ der Dimension
$s\ge5$ eine Ungleichung der Form $$ \abs{P_\E(x)}\le C_s\,
B(r_{\max},r_{\min})\,x^{s/2-1}\,,\eqno(1.2)  $$ wobei $C_s$ nur von
$s$ abh\"{a}ngt und $B(r_{\max},r_{\min})$ mittels des minimalen und
maximalen Hauptkr\"ummungsradius des Ellipsoids explizit angegeben
wird. \ssk F\"{u}r die Einheitskreisscheibe wurde von E.~Kr\"{a}tzel
[9], Satz 5.12, die Absch\"atzung
$$ \abs{P(x)} \le 38\,x^{1/3}+704\,x^{1/4} + 11 \eqno(1.3) $$
erzielt. Allgemeiner konnte er in [10] f\"{u}r einen ebenen Bereich
$\B$ mit glattem Rand von stetiger, beschr\"{a}nkter, nicht
verschwindender Kr\"{u}mmung
$$ \abs{P_\B(x)} \le 48\,\(r_{\max}^2\,x\)^{1/3} + \(703\,\sqrt{r_{\max}}
+{3\,r_{\max}\over5\,\sqrt{r_{\min}}}\)\,x^{1/4} + 11 \eqno(1.4)
$$ beweisen, wobei $r_{\max}$, $r_{{\min}}$ die Extremwerte des
Kr\"{u}mmungsradius der Randkurve bezeichnen. Gemeinsam untersuchten
die Autoren in [11] den Gitterrest von Ellipsenscheiben und
Ellipsoiden im $\R^3$ der Gestalt $Q(\b u)\le x$, \ $Q$ eine positiv
definite quadrati\-sche Form mit unimodularer
Koeffizientendeterminante. Die erzielten Schranken lauten $$
\abs{P_\E(x)} \le 8.46\,x^{1/3} + \hbox{ Terme kleinerer Ordnung in
} x \eqno(1.5)$$ in der Ebene und
$$ \abs{P_\E(x)} \le 14\,x^{3/4} + \hbox{ Terme kleinerer Ordnung
in } x \eqno(1.6)$$ im $\R^3$. Dieselbe Arbeit [11] behandelt auch
den Fall eines $\R^3$-Rotationsk\"{o}rpers (bez\"{u}g\-lich einer
Koordinatenachse) mit glattem Rand von durchwegs beschr\"{a}nkter,
nicht verschwindender Gau{\ss}scher Kr\"{u}mmung, mit einem Ergebnis
der Form
$$ \abs{P_\K(x)}\le C(r_{\max},r_{\min})\,x^{3/4}+ \hbox{ Terme
kleinerer Ordnung in } x\,. \eqno(1.7)$$ Dabei ist
$C(r_{\max},r_{\min})$ ein expliziter Ausdruck, der von der
erzeugenden Kurve abh\"{a}ngt, insbesondere von den Extremwerten
ihres Kr\"{u}mmungsradius. F\"{u}r den Spezialfall der Kugel ist
$C(r_{\max},r_{\min})=86$. \bsk\bsk

{\bf 2.~Gegenstand und Resultat dieser Arbeit. } Wir betrachten hier
den speziellen Fall eines bez\"{u}glich der dritten Koordinatenachse
rotationssymmetrischen Ellipsoids $\E$ in Mittelpunktslage. Passend
normiert (so dass das Volumen $\vol$ betr\"{a}gt), l\"{a}sst dieses
sich in der Gestalt
$$ \E:\quad Q(\b u) := {u_1^2+u_2^2\over a}+a^2\,u_3^2\ \le\,1 \eqno(2.1)
$$ darstellen, mit einer Konstanten $a>0$. F\"{u}r gro{\ss}es reelles $x$
werden wir f\"{u}r die Gitter-Diskrepanz des linear
vergr\"{o}{\ss}erten Ellipsoids $Q(\b u)\le x$ eine Absch\"{a}tzung
der Form $$ \abs{P_\E(x) } \le 1237\,a^{1/8}\,x^{11/16}\,\L^{3/8}\ +
\hbox{ Terme kleinerer Ordnung in } x \eqno(2.2) $$ beweisen, wobei
f\"{u}r den logarithmischen Faktor stets die Abk\"{u}rzung $$ \L :=
\log(100x)+|\log a| \eqno(2.3)
$$ ben\"{u}tzt wird. Diese Schranke ist bez\"{u}glich des Exponenten von $x$
deutlich sch\"{a}rfer als (1.6) und (1.7) (${11\over16}=0.6875$
gegen\"{u}ber ${3\over4}=0.75$) und unterschreitet damit den
Landau-Hlawka'schen (vgl.~[6]) Fehler-Exponenten $s(s-1)\over2(s+1)$
f\"{u}r Dimension $s$. \ssk Unsere Methode st\"{u}tzt sich auf jene
von F.~Chamizo [2], der f\"{u}r einen allgemeinen
$u_3$-rotationssymmetrischen K\"{o}rper $\K$ des $\R^3$ mit glattem
Rand von beschr\"{a}nkter, nirgends verschwindender Gau{\ss}scher
Kr\"{u}mmung $$ P_\K(x) = O_{\K,\ep}\(x^{11/16+\ep}\)\qquad(\ep>0)
\eqno(2.4) $$ bewies. Allerdings ben\"{o}tigt Chamizo eine
technische Bedingung an die dritte Ableitung der erzeugenden
Funktion der rotierenden Kurve, die kaum geometrisch interpretierbar
ist. Au{\ss}erdem wird in seinem Argument ein wesentlicher Fall (bei
uns mit $\II$ bezeichnet, vgl.~ab Formel (4.21) unten) nicht
ausgef\"{u}hrt. Die Absch\"{a}tzung (2.4) wurde \"{u}brigens f\"{u}r
den Spezialfall der Kugel bereits von I.M.~Vinogradov [15] erzielt.
Er versch\"{a}rfte sein Ergebnis sp\"{a}ter bis zu $O\(x^{2/3}(\log
x)^6\)$ [16]. Die beste bekannte Schranke f\"{u}r die Kugel stammt
heute von D.R.~Heath-Brown [5] (nach Vorarbeit F.~Chamizo und
H.~Iwaniec [3]) und lautet $$ P(x) = \O{x^{21/32+\ep}}\,.\eqno(2.5)
$$

\ssk Selbst eine nicht-effektive Fassung von (2.2) mit
$O_a\(x^{11/16}(\log x)^{3/8}\)$ w\"{u}rde f\"{u}r das
Rotationsellipsoid alles Bekannte versch\"{a}rfen. Unser Zugang ist
optimiert f\"{u}r den Exponenten von $x$, also am effizientesten
f\"{u}r den Fall, dass $|\log a|$ klein gegen $\log x$ ist. In
unserem Beweis werden zwei von $x$ und $a$ abh\"{a}ngige Parameter
$$ y:=73.6\,a^{1/8}\,x^{3/16}\,\L^{3/8}\,,\quad z:= {0.3852\over y}
\,\sqrt{x+2y} \eqno(2.6) $$ auftreten. Wir ben\"{o}tigen die
technische Voraussetzung, dass $$ 1\le y\le{x\over3} \eqno(2.7)
$$ gilt. Diese ist sicher erf\"{u}llt z.B.~f\"{u}r $$ {1\over x}\le a\le
x\,,\quad x\ge15\,000\,. \eqno(2.8) $$ Unser Resultat lautet wie
folgt. \bsk

\vbox{{\bf Satz. } {\it F\"{u}r positive reelle Zahlen $a, x$,
f\"{u}r welche die Bedingung $(2.7)$ erf\"{u}llt ist, gen\"{u}gt die
Gitter-Diskrepanz $P_\E(x)$ des Rotationsellipsoids $$ Q(\b u) =
{u_1^2+u_2^2\over a} +a^2\,u_3^2\le x
$$ der Ungleichung $$ \eq{\abs{P_\E(x)}\le &\ 1237\,a^{1/8}\,x^{11/16}\,\L^{3/8} +
12\,a^{-69/64}\,x^{81/128}\,\L^{145/64} \cr & +
\(a^{1/4}(134\,\L^{5/2}+543\,\L^{1/4})+20\,{\L^{5/2}\over
a^{1/4}}\)\,x^{5/8} + 12\,a^{-39/64}\,x^{75/128}\,\L^{139/64}\cr & +
\(268\,{\L\over
a}+159\,\al_0^3+2000\)\,x^{1/2}+\al_0^3\,(4.4\,\L+104.5)\,, \cr}
$$ wobei $\L$ durch $(2.3)$ erkl\"{a}rt ist und $\al_0:=\max(a,1/\sqrt{a})$
gilt.}}

\bsk \msk

{\bf3.~Einige Hilfss\"{a}tze.} Die folgenden Ungleichungen sind
qualitativ wohlbekannt. Es geht uns hier um Versionen mit expliziten
Konstanten. \bsk

{\bf Lemma 1.} \quad {\it Wie \"{u}blich bezeichne $r(n)$ die Anzahl
der Darstellungen von $n$ als Summe der Quadrate zweier ganzer
Zahlen.\ssk {\bf(a) } F\"{u}r $x\ge1$ gilt dann
$$ R_1(x):=\sum_{1\le n\le x} r(n) \le \ 4x $$ und $$ R_{1,2}(x):=
\sum_{x<n\le2x} r(n) \le \ 4.8\,x\,. $$ {\bf(b)} F\"{u}r alle
$n\in\N^*$ ist $$ r^2(n) \le \sum_{mk=n}r(m)r(k)\,. $$ {\bf(c)}
F\"{u}r $x\ge1$ gilt $$ R_2(x):=\sum_{x<n\le2x}r^2(n)\le \
19.2\,x\log(2e^2\,x)\,.
$$ } \bsk

{\bf Beweis. } Die \"{u}bliche elementare Idee, jedem Gitterpunkt
$\b m\in\Z^2$ ein Einheitsquadrat mit Mittelpunkt $\b m$ zuzuordnen,
ergibt
$$R_1(x)\le \pi\(\sqrt{x}+{1\over\sqrt{2}}\)^2-1\,. $$
Dies ist f\"{u}r $x\ge29$ tats\"{a}chlich $\le4x$. Man verifiziert
leicht (z.B.~mit Hilfe von {\it Derive } [13]), dass  $$ \max_{1\le
x\le29}{R_1(x)\over x}=\max_{n=1,\dots,29} {R_1(n)\over n} = 4
 $$ gilt. Ebenso folgt $$ R_{1,2}(x)
 \le\pi\(\(\sqrt{2x}+{1\over\sqrt{2}}\)^2-\(\sqrt{x}-{1\over\sqrt{2}}\)^2\)\,,
 $$ dies ist $\le 4.8\,x$ f\"{u}r $x\ge42$. Direkte Berechnung zeigt
 $$\max_{n=2,\dots,90}{R_{1,2}(n/2)\over n/2}=4.8\,. \eqno{\qed}$$
 {\bf(b)} ist Lemma E, Formel (2.1), in Nowak [12]. \qed \ssk {\bf(c)} Es folgt $$
R_2(x)\le \sum_{x< mk\le2x}r(m) r(k) \le 2 \sum_{1\le
m\le\sqrt{2x}}r(m)\,R_{1,2}\({x\over m}\)\le $$ $$ \le
9.6\,x\sum_{1\le m\le\sqrt{2x}}{r(m)\over m} = 9.6\,x
\Int_{1-}^{\sqrt{2x}}{1\over w}\,\d R_1(w)=
9.6\,x\,{R_1(\sqrt{2x})\over\sqrt{2x}}+
9.6\,x\Int_1^{\sqrt{2x}}{R_1(w)\over w^2}\,\d w \le$$ $$ \le
38.4\,x+38.4\,x\Int_1^{\sqrt{2x}}{\d w\over w} =
38.4\,x\(1+\h\log(2x)\)=19.2\,x \log(2e^2\,x)\,. \eqno{\qed} $$
\bsk\msk

{\bf Lemma 2. } {\it F\"{u}r Elemente des $\R^3$ bezeichne
$\nm\cdot$ die "Zylindernorm"
$\nm{(w_1,w_2,w_3)}:=\max\(\sqrt{w_1^2+w_2^2},|w_3|\)$. Dann gilt:
\ssk \b{(a) } F\"{u}r alle $x>0$ ist $$R_3^*(x) := \#\{\b m\in\Z^3:\
0<\nm{\b m }^2\le x\} \le 14\,x^{3/2}\,.$$ \ssk \b{(b) } F\"{u}r $\b
m=(m_1,m_2,m_3)\in\Z^3$, $a>0$, sei $g_*(\b
m):=\max(\sqrt{a(m_1^2+m_2^2)},|m_3|/a)$  und $g_0:= 1/\al_0
=\min(\sqrt{a},1/a)$ definiert. Dann folgt f\"{u}r $Z>0$
$$ \sum_{0<g_*({\b m})\le Z}(g_*({\b m}))^{-3} \le {42\over g_0^3}\,
\log_+(1.4\,Z/g_0)\,,
$$ wobei $\log_+ := \max(\log,0)$ bedeutet. \ssk
\b{(c) } Ebenso gilt f\"{u}r $Z>0$, $\al>3$, $$ \sum_{g_*({\b m})>
Z} (g_*({\b m})^{-\al} \le{14\al\over(\al-3)g_0^3}\,Z^{3-\al}\,. $$
} \bsk {\bf Beweis. } \b{(a) } Adaption des \"{u}blichen Arguments
ergibt
$$ R_3^*(x) \le
\pi\(\sqrt{x}+{1\over\sqrt{2}}\)^2\(2\sqrt{x}+1\)-1\,. $$ Dies ist
$\le14\,x^{3/2}$ f\"{u}r $x\ge5$. Die Werte $x= 1, 2, 3, 4$ werden
einzeln \"{u}berpr\"{u}ft. \qed \ssk \b{(b) } Es ist $g_*(\b m)\ge
g_0 \nm{\b m}$, daher gilt $$ G(u):=\sum_{0<g_*(\b m)\le u}1\le
R_3^*(u^2/g_0^2)\ \cases{\le{14\,u^3 / g_0^3}& immer, \cr = 0 & wenn
$u<g_0$.\cr} $$ Durch partielle Summation folgt, falls $Z\ge g_0$,

\vbox{$$ \sum_{0<g_*({\b m})\le Z} (g_*({\b m}))^{-3} =
\Int_{0+}^{Z} u^{-3}\d G(u) = Z^{-3}\,G(Z)+ 3 \Int_{0+}^{Z}
u^{-4}G(u)\d u $$ $$  \le {14\over g_0^3}\(1+3\Int_{g_0}^Z{\d u\over
u }\)= {42\over g_0^3}\,\log(e^{1/3}\,Z/g_0)\,.$$}

F\"{u}r $Z<g_0$ ist diese Summe nat\"{u}rlich 0. \qed \ssk

\b{(c) } Ganz analog erh\"{a}lt man
$$\sum_{g_*({\b m})> Z} (g_*({\b m}))^{-\al} = \Int_{Z+}^\infty
u^{-\al}\d G(u)$$ $$ \le \al \Int_Z^{\infty} u^{-\al-1}G(u)\d u \le
{14\al\over g_0^3}\Int_Z^\infty u^{2-\al}\d
u={14\al\over(\al-3)g_0^3}\,Z^{3-\al}\,.\eqno{\qed}$$

\bsk  \msk

{\bf 4.~Beweis des Satzes.\msk 4.1.~Vorbereitung der
Absch\"{a}tzung. } Wir verwenden die auf E.~Landau
zur\"{u}ck\-gehende Mittelwert-Technik, die E.~Hlawka f\"{u}r
allgemeine konvexe K\"{o}rper in [6] im Detail ausgef\"{u}hrt hat.
Dazu definieren wir f\"{u}r jede Funktion $F$, die auf jedem
kompakten Teilintervall von $[0,\infty[$ st\"{u}ckweise stetig und
beschr\"{a}nkt ist,
$$ F_\z(t):= \Int_0^t \(\Int_0^{t_1}F(t_2)\d t_2\)\d t_1 \qquad(t>0)
\eqno(4.1)
$$ und, f\"{u}r $x>0$, $|u|\le{1\over3}x$,
$$ D_{x,u}^\z(F) := F_\z(x+2u)-2F_\z(x+u)+F_\z(x) =
\Int_x^{x+u} \(\Int_{t_1}^{t_1+u}F(t_2)\d t_2\)\d t_1\,. \eqno(4.2)
$$ Es bezeichne $V(t)=\vol t^{3/2}$ das Volumen des Ellipsoids
$Q \le t$. F\"{u}r $0<y\le{1\over3}x$ folgt wegen
$$ A(x)\le y^{-2} D_{x,y}^\z (A)\,,\quad A(x)\ge y^{-2}
D_{x,-y}^\z (A)$$ unmittelbar $$ |P(x)|\le \max_\pm \abs{y^{-2}
D_{x,\pm y}^\z(V)-V(x)} + \max_\pm\abs{y^{-2} D_{x,\pm y}^\z(P)}\,.
\eqno(4.3) $$ Nun ist $$D_{x,\pm y}^\z(V)=\vol\,
{4x^{7/2}\over35}\,\phi\({\pm y\over x}\)\,,\quad
\phi(w):=(1+2w)^{7/2}-2(1+w)^{7/2}+1\,. $$ Nach Taylor folgt f\"{u}r
$|w|\le{1\over3}$ $$ \abs{\phi(w)-{35\over4}\,w^2}\le{|w|^3\over6}\,
\max_{|v|\le1/3}\abs{\phi'''(v)} \,,$$ daher $$ \abs{y^{-2} D_{x,\pm
y}^\z(V)-V(x)}= \vol\, {4x^{7/2}\over35\,y^2}\,\abs{\phi\({\pm
y\over x}\)- {35\over4}\,\({\pm y\over x}\)^2}$$ $$ \le
{8\pi\over315}\,{x^{7/2}\over y^2}\,{y^3\over
x^3}\,\max_{|v|\le1/3}\abs{\phi'''(v)} \le 8.4\,x^{1/2}y\,,$$ also
mit (4.3) $$ |P(x)| \le 8.4\,x^{1/2}y + \max_\pm\abs{y^{-2} D_{x,\pm
y}^\z(P)}\,. \eqno(4.4) $$ Nun l\"{a}sst sich $D_{x,\pm y}^\z(P)$
mittels der Poisson'schen Formel durch eine absolut konvergente
Reihe darstellen (vgl.~z.B.~Hlawka [6]). Es ist
$$ D_{x,\pm y}^\z(P) = \sum_{\b o\ne\b m\in\Zi^3} D_{x,\pm y}^\z
\(I(\b m, \cdot)\)\,,\quad\quad I(\b m,t):= \Int_{Q(\b u)\le t} e(\b
m\,\b u)\d\b u\,, $$ mit $e(w)=e^{2\pi i w}$ wie \"{u}blich. Diese
Reihe wird nach der Gr\"{o}{\ss}e von $g_*(\b
m):=\max(\sqrt{a(m_1^2+m_2^2)}, |m_3|/a)$ in zwei Teile zerlegt:
$$ D_{x,\pm y}^\z(P) = \sum_{0< g_*(\b m)\le z} D_{x,\pm y}^\z
\(I(\b m, \cdot)\) + \sum_{g_*(\b m) > z} D_{x,\pm y}^\z \(I(\b m,
\cdot)\) =: S_{\rm I} + S_{\rm II}\,, $$ wobei $z>0$ noch
verf\"{u}gbar bleibt. Zur Absch\"{a}tzung von $S_{\rm I}$ wird die
zweite der in (4.2) enthaltenen Darstellungen verwendet, f\"{u}r
$S_{\rm II}$ hingegen die erste. So erh\"{a}lt man $$
\eq{\abs{S_{\rm I}} &\le y^2 \max_{|t-x|\le2y}\abs{\sum_{0<g_*(\b m)
\le z}I(\b m,t)}\,,\cr \abs{S_{\rm II}} &\le 4
\max_{|t-x|\le2y}\abs{\sum_{g_*(\b m)> z}I_\z(\b m,t)}\,.\cr  } $$
Die einzelnen Summanden lassen sich explizit auswerten. Es ergibt
sich
$$ \eq{I(\b m,t) &= t^{3/2} \Int_{Q(\b u)\le1}e(\b m \sqrt{t}\,
\b u)\d\b u = t^{3/2} \Int_{|\b v|_2\le1}e(g(\sqrt{t}\,\b m)v_1)\d\b
v =\cr &= t^{3/2}\Int_{-1}^1 e(g(\sqrt{t}\,\b
m)v_1)\(\Int_{v_2^2+v_3^2\le1-v_1^2} \d(v_2,v_3)\)\d v_1=\cr &=
 - {\sqrt{t}\,\cos(2\pi g(\b m)\sqrt{t})\over\pi g^2(\b m)} +
 {\sin(2\pi g(\b m)\sqrt{t})\over2\pi^2 g^3(\b m)}\,,\cr }
\eqno(4.5)  $$ wobei $g=\sqrt{Q^{-1}}$ ist, $Q^{-1}$ die inverse
Form zu $Q$, also $g(\b m)=\sqrt{a(m_1^2+m_2^2)+m_3^2/a^2}$. Durch
zweimalige Integration folgt weiter
$$ \eq{I_\z(\b m,t) &= {t^{3/2}\cos(2\pi g(\b m)\sqrt{t})\over\pi^3 g^4(\b
m)} - {3t\,\sin(2\pi g(\b m)\sqrt{t})\over\pi^4 g^5(\b m)}
 \cr &- {15\sqrt{t}\,\cos(2\pi g(\b m)\sqrt{t})\over4\pi^5 g^6(\b m)}
 + {15\sin(2\pi g(\b m)\sqrt{t})\over8\pi^6 g^7(\b m)}\,.  \cr } \eqno(4.6)
 $$ Wir werden von (4.5) und (4.6) nur jeweils die ersten Terme
 genau behandeln, die \"{u}brigen trivial absch\"{a}tzen. Wir ben\"{u}tzen
 Lemma 2 und erhalten wegen $g\ge g_*$ $$
 {1\over2\pi^2}\sum_{0<g_*(\b m)\le z} g^{-3}(\b m)\le {21\over\pi^2
 g_0^3}\,\log_+(1.4\,z/g_0) $$ und, f\"{u}r $t\ge1$,
 $$ {4\over y^2}\sum_{g_*(\b m)>z}\({3t\over\pi^4} g^{-5}(\b m)+
 {15\sqrt{t}\over4\pi^5}g^{-6}(\b m)+{15\over8\pi^6}
 g^{-7}(\b m)\)\le {6t\over g_0^3 y^2 z^2} \,.$$ F\"{u}r $g_*(\b
 m)>x/y^2$ sch\"{a}tzen wir auch den Hauptteil von (4.6) trivial ab:
 $$ {4t^{3/2}\over\pi^3 y^2} \sum_{g_*(\b m)>x/y^2} g^{-4}(\b
 m)\le {224\,t^{3/2}\over\pi^3 g_0^3\, x}\,,  $$ wieder nach
 Lemma 2. Als n\"{a}chstes werden noch die $\b
 m\in\Z^3$ mit $m_1^2+m_2^2\le20$ und jene mit $|m_3|\le40$ extra
 behandelt, und zwar durchwegs mittels (4.5), die Unterscheidung
 nach der Gr\"{o}{\ss}e von $g_*(\b m)$ erfolgt hier nicht. Wir erhalten,
 unter Verwendung von Lemma 1 (a),
 $$ \sum_{\b m:\ m_1^2+m_2^2\le20} \({\sqrt{t}\over\pi}\,g^{-2}(\b m)+
 {1\over2\pi^2}\,g^{-3}(\b m)\) $$ $$ \eq{&\le81\sum_{m_3\ne0} \({\sqrt{t}\over\pi}\,{a^2\over m_3^2}+
 {1\over2\pi^2}\,{a^3\over|m_3|^3}\) + 81\,\({\sqrt{t}\over\pi}\,{1\over a}+
 {1\over2\pi^2}\,{1\over a^{3/2}}\)\cr &\le {81\over\pi}\,{\sqrt{t}\over g_0^2}(2\zeta(2)+1)+
{81\over2\pi^2}\,{1\over g_0^3}(2\zeta(3)+1)\le 111\,{\sqrt{t}\over
g_0^2}+{14\over g_0^3}\,. \cr }$$ \"{A}hnlich ergibt sich $$
\sum_{0<g_*(\b m)\le x/y^2:\atop |m_3|\le40, (m_1,m_2)\ne(0,0)}
\({\sqrt{t}\over\pi}\,g^{-2}(\b m)+
 {1\over2\pi^2}\,g^{-3}(\b m)\) $$ $$ \le 81 \sum_{0<a(m_1^2+m_2^2)\le
 x^2/y^4} \({\sqrt{t}\over\pi}\,{1\over a(m_1^2+m_2^2)}+
 {1\over2\pi^2}\,{1\over a^{3/2}(m_1^2+m_2^2)^{3/2}}\)$$ $$ \le
 {324\,\sqrt{t}\over\pi\,a}\,\log(e\,x^2/(ay^4))+{486\over\pi^2\,a^{3/2}}\,,
 $$ da wegen Lemma 1 (a) mit partieller Summation $$ \sum_{0<n\le Z}{r(n)\over n}\le 4 \log(e\,Z)\,,\qquad
\sum_{n>0}{r(n)\over n^{3/2}}\le12  $$ gilt.
 Im Folgenden k\"{o}nnen wir uns also auf die $\b m\in\Z^3$ mit
$$ m_1^2+m_2^2>20\,,\quad |m_3|>40 \eqno(4.7) $$ beschr\"{a}nken, was in den betreffenden Summen als
$\disp\sum^{(4.7)}$ angedeutet wird. \ssk Kombination aller dieser
Schranken mit (4.4) - (4.6) ergibt
 $$ \eq{|P(x)| \le &\ \max_{|t-x|\le2y}\abs{{\sqrt{t}\over\pi}
 \sum_{0<g({\b m})\le z}^{(4.7)} {e(g(\b m)\sqrt{t})\over g^{2}(\b m)}} +
 \max_{|t-x|\le2y}\abs{{4\,t^{3/2}\over\pi^3
 y^2} \sum_{z<g({\b m})\le x/y^2}^{(4.7)} {e(g(\b m)\sqrt{t})\over g^{4}(\b m)}} \cr & +
 8.4\,x^{1/2}y + {2.2\over g_0^3}\,\log_+(1.4\,z/g_0) + {6t^*\over g_0^3 y^2 z^2}\cr & +
 {224\,(t^*)^{3/2}\over\pi^3 g_0^3\, x}+111\,{\sqrt{t^*}\over g_0^2}+{14\over g_0^3}+{207\over a}\,
 \sqrt{t^*}\,\L+{50\over a^{3/2}}\,,\cr }  \eqno(4.8) $$
 wobei $t^*:=x+2y$ gesetzt wurde.

 \bsk

 {\bf4.2.~Absch\"{a}tzung der Exponentialsummen. } Es sei
 $f(n,m):=\sqrt{an+m^2/a^2}=g(\sqrt{n},0,m)$, dann
 legt (4.8) nahe, Teilsummen der Gestalt $$
 S^{[j]}(N,M):=\sum_{N<n\le2N}r(n) \sum_{M<m\le\wz M} (f(n,m))^{-j}\, e(f(n,m)\sqrt{t})
 \eqno(4.9)$$ f\"{u}r reelle $N,M\ge1$ und $j=2$ oder $j=4$ zu betrachten. F\"{u}r $u\in[N,2
 N]$, $w\in[M,\wz M]$, sei $$  E_{N,M}(u,w):=
 \sum_{N<n\le u}r(n) \sum_{M<m\le w}
 e(f(n,m)\sqrt{t})\,,\eqno(4.10)
 $$ dann folgt durch zweimalige partielle Summation $$
 \eq{& S^{[j]}(N,M)= \ f^{-j}(2 N,\wz M)E_{N,M}(2 N,\wz M) \cr -& \Int_N^{2
 N} {\partial\over\partial u}\(f^{-j}(u,\wz M)\) E_{N,M}(u,\wz M)\d u \cr -
& \Int_M^{\wz M} {\partial\over\partial w}\(f^{-j}(2 N,w)\)E_{N,M}(2
N,w)\d
 w \cr + \Int_M^{\wz M}&\(\Int_N^{2 N}  {\partial^2\over\partial u\partial
 w}\(f^{-j}(u,w)\)E_{N,M}(u,w) \d u\)\d w\,.\cr  }  $$ (Man
 vgl.~Kr\"{a}tzel [8], Theorem 1.6, f\"{u}r einen entsprechenden allgemeinen
 Satz.) Da die partiellen Ableitungen von $f^{-j}(u,w)$ nicht verschwinden,
 folgt daraus $$ \abs{S^{[j]}(N,M)}\le {4\over f^j(N,M)}\,\sup_{u\in[N,2
 N]\atop w\in[M,\wz M]}\abs{E_{N,M}(u,w)}\,.\eqno(4.11) $$ Im
 Folgenden seien $U\in[N,2N]$, $W\in[M,\wz M]$ ganze Zahlen, f\"{u}r die
 das Supremum in (4.11) angenommen wird. Die Exponentialsumme
 $E=E_{N,M}(U,W)$ wird durch einen sog.~Weylschen Schritt weiter
 abgesch\"{a}tzt. Es sei $H$ eine ganze Zahl mit $$ 10\le H\le\h M \,.
 \eqno(4.12) $$ Dann gilt $$
 H\,E = \sum_\m \ \sum_\n r(n) \sum_{k\in[1,H]:\ \mk}
 e\(f(n,m+k)\wt\)\,, $$ daher ergibt zweimalige Anwendung der
 Cauchyschen Ungleichung
 $$ H^2|E|^2\le M \sum_\m \abs{\sum_\n r(n) \sum_{k\in[1,H]:\ \mk}
 e\(f(n,m+k)\wt\)}^2 $$ $$ \le M\sum_\m \(\sum_\n r(n)\abs{\sum_{k\in[1,H]:\ \mk}
 e\(f(n,m+k)\wt\)}\)^2 $$ $$ \le M\(\sum_\n r^2(n)\)\sum_\m \ \sum_\n \abs{\sum_{k\in[1,H]:\ \mk}
 e\(f(n,m+k)\wt\)}^2 $$ $$ \le M R_2(N)\sum_\m \sum_\n \sum_{k_1,k_2\in[1,H]:\atop m+k_1, m+k_2\in]M,W]}
 e\((f(n,m+k_1)-f(n,m+k_2))\wt\) $$
 Die Diagonalterme ergeben einen  Beitrag $\le\h M^2 H N R_2(N)$, in der verbleibenden Summe setzen
 wir $m'=m+\min(k_1,k_2)$ und $h=|k_1-k_2|$ als neue Summationsvariable und
 erhalten:
 $$ \eq{& H^2|E|^2 \le \h M^2 H N R_2(N) \cr
+ 2M H R_2(N)&\sum_{h\in[1,H], m'\in]M,W]:\atop m'+h\in]M,W]}
 \abs{\sum_\n e\((f(n,m'+h)-f(n,m'))\wt\)}\,.\cr}  $$

\vbox{Lemma 1 (c) vereinfacht dies zu
 $$ \eq{H|E|^2 &\le 9.4\, M^2 N^2 \log(2e^2 N)
+ 38.4\, M N\log(2e^2N)\times\cr \times&\sum_{h\in[1,H],
m\in]M,W]:\atop m+h\in]M,W]}
 \abs{\sum_\n e\((f(n,m+h)-f(n,m))\wt\)}\,,\cr}  \eqno(4.13)$$
 wenn wir statt $m'$ wieder $m$ schreiben.} \ssk
 Es verbleibt, f\"{u}r $F(\tau)=F_{h,m;t}(\tau):=(f(n,m+h)-f(n,m))\wt$ die
 innere Exponentialsumme $\disp\sum_{N<n\le U}e(F(n))$
 abzusch\"{a}tzen, wobei $h\in[1,H]$, $m, m+h \in ]M,W]$ voraus\-gesetzt wird.
 Wir verwenden dazu eine effektive Version der einfachsten Van der
 Corput'schen Schranke: vgl.~Kr\"{a}tzel [9], S.~16, Formel (1.10). Aus
 $$
 F(\tau)=F_{h,m;t}(\tau)=\(\sqrt{a\tau+(m+h)^2/a^2}-\sqrt{a\tau+m^2/a^2}\)\wt
 $$ folgt $$ F''(\tau)={3\over4}\wt \Int_m^{m+h}
 {\xi\,\d\xi\over(a\tau+\xi^2/a^2)^{5/2}}\,, $$ daher ist $$ \min_{N\le\tau\le
 U} F''(\tau)\ge
 {3\over4}{hm\wt\over(a\tau+(m+h)^2/a^2)^{5/2}}\ge
 {3\over4}{hM\wt\over(2aN+2M^2/a^2)^{5/2}}:=\Lambda $$ und $$
 \max_{N\le\tau\le U} F''(\tau)\le
 {3\over4}\wt\,{h(m+h)\over(a\tau+m^2/a^2)^{5/2}} \le {3\over4}\wt\,
 {h\wz\,M\over(aN+M^2/a^2)^{5/2}}=8\Lambda\,. $$  Daher impliziert
 Formel (1.10) in [9], S.~16,
 $$  \eq{&\abs{\sum_{N<n\le
 U}e(F_{h,m;t}(n))}\le40(U-N)\sqrt{\Lambda}+{11\over\sqrt{\Lambda}}\cr
 &\le {\sqrt{3}\over\root4\of2}\,{10N\,\sqrt{hM}\,\root4\of
 t\over(f(N,M))^{5/2}} + {\root4\of2\over\sqrt{3}}\,{44\over\sqrt{hM}\,
 \root4\of t}\,(f(N,M))^{5/2}\,. \cr  } \eqno(4.14) $$  Wir
 verwenden dies in (4.13): Die Summation \"{u}ber $m$ ergibt einen Faktor
 $\le\h M$, jene \"{u}ber $h$ Faktoren $\le H^{3/2}$ bzw.~$\le2\sqrt H$. So vereinfacht sich (4.13) zu
 $$ \eq{H\,|E|^2\le&\ 9.4\,M^2 N^2 \log(15N) + {280\,M^{5/2}\,N^2\log(15N)\,H^{3/2}\root4\of t
\over(f(N,M))^{5/2}}\cr & + 1161\,M^{3/2}N \log(15N)\,\sqrt{H}\,
(f(N,M))^{5/2}\,t^{-1/4}\,.\cr }  $$ Division durch $H\,\log(15\,N)$
und Verwendung von $\sqrt{A+B}\le\sqrt{A}+\sqrt{B}$ ergibt
$$ {|E|\over\sqrt{\log(15N)}} \le
3.1\,{MN\over\sqrt{H}}+16.8\,{M^{5/4}\,N\,H^{1/4}t^{1/8}
\over(f(N,M))^{5/4}}+ 34.1\,{M^{3/4}N^{1/2}\, (f(N,M))^{5/4}\over
H^{1/4}\,t^{1/8}}\,. $$ Wir balanzieren hier die ersten beiden Terme
gegeneinander aus und erhalten $$ H =
\left[0.1\,M^{-1/3}(f(N,M))^{5/3}\,t^{-1/6}\right]\,. \eqno(4.15)$$
Wir betrachten nun zun\"{a}chst den Fall, dass dieser Wert f\"{u}r
$H$ die Bedingungen (4.12) erf\"{u}llt. Diese Einschr\"{a}nkung wird
in den auftretenden Gr\"{o}{\ss}en durch das Superskript $\I$
symbolisiert. \ssk Da f\"{u}r $\al\ge10$ sicher $[\al]\ge0.9\,\al$
gilt, folgt aus obiger Absch\"{a}tzung, wenn wir kurz $f:=f(N,M)$
schreiben,
$$ {\abs{E_{N,M}(U,W)}\over\sqrt{\log(15N)}} \le
20\,M^{7/6}N\,f^{-5/6}\,t^{1/12}+ 63\,M^{5/6} N^{1/2}
f^{5/6}\,t^{-1/12}\,.\eqno(4.16)  $$ Nach (4.11) erhalten wir
f\"{u}r $j=2$ oder $j=4$ $$
{\abs{S^{[j]}(N,M)}\over\sqrt{\log(15N)}}\le
80\,M^{7/6}N\,f^{-5/6-j}\,t^{1/12}+ 252\,M^{5/6} N^{1/2}
f^{5/6-j}\,t^{-1/12}\,.\eqno(4.17)  $$ Zur Behandlung der ersten
Summe in (4.8) verwenden wir dies mit $j=2$ sowie
$N=N_r=2^{-r}\,z^2/a$, $M=M_s=2^{-s/2}\,a z$ mit $r,s\ge0$ und ganz.
Dann ist $f(N_r,M_s)=(2^{-r}+2^{-s})^{1/2}\,z$. Setzen wir noch
$\L:=\log(100\,x)+|\log a|$, dann folgt
$$ \eq{\abs{\sum_{0<g_*({\b m})\le z}\I {e(g(\b m)\sqrt{t})\over g^{2}(\b
m)}}\,\L^{-1/2}\le & 160\,a^{1/6}\,t^{1/12}\,z^{1/3}\sum_{r,s\,
\ge0}2^{-r-7s/12}(2^{-r}+2^{-s})^{-17/12}\cr + &
504\,a^{1/3}\,t^{-1/12}\,z^{2/3}\sum_{r,s\,
\ge0}2^{-r/2-5s/12}(2^{-r}+2^{-s})^{-7/12}\,,\cr}
$$ wenn man ber\"{u}cksichtigt, dass jedem Intervall $]M, \wz M]$ ein Intervall
$[-\wz M,-M[$ entspricht. Wegen
$\root3\of{\al(\be/2)^2}\le(\al+\be)/3$ mit $\al=2^{-s}, \be=2^{-r}$
ist die Konvergenz der ersten Reihe evident, analog wegen
$\sqrt{\al\be}\le(\al+\be)/2$ die der zweiten. Numerische Berechnung
mit einem Computeralgebraprogramm zeigt
$$\sum_{r,s=0}^\infty{2^{-r-7s/12}\over(2^{-r}+2^{-s})^{17/12}}<23.8\,,\qquad
\sum_{r,s=0}^\infty{2^{-r/2-5s/12}\over(2^{-r}+2^{-s})^{7/12}}<27\,,
$$ damit erhalten wir $$ \abs{\sum_{0<g_*({\b m})\le z}\I {e(g(\b m)\sqrt{t})\over g^{2}(\b
m)}}\le\L^{1/2}\(3808\,a^{1/6}\,t^{1/12}\,z^{1/3}+13\,608\,a^{1/3}\,t^{-1/12}\,z^{2/3}\)\,.
\eqno(4.18)$$ F\"{u}r die zweite Summe in (4.8) ben\"{u}tzen wir
(4.17) mit $j=4$. $r, s$ durchlaufen jetzt endliche Mengen ganzer
Zahlen, wobei aus $g_*(\b m)>z$ die Bedingung $\min(r,s)\le0$ folgt.
Wir erhalten
$$ \eq{\abs{\sum_{z<g_*({\b m})\le x/y^2}\I {e(g(\b
m)\sqrt{t})\over g^{4}(\b m)}}\,\L^{-1/2}\le &
160\,a^{1/6}\,t^{1/12}\,z^{-5/3}\sum_{r,s}{2^{-r-7s/12}\over(2^{-r}+2^{-s})^{29/12}}\cr+
&504\,a^{1/3}\,t^{-1/12}\,z^{-4/3}\sum_{r,s}{2^{-r/2-5s/12}\over(2^{-r}+2^{-s})^{19/12}}\,.\cr}
$$ Die Summen \"{u}ber $r,s$ werden nun aufgeteilt in der Form $$
\sum_{\min(r,s)\le0} = \sum_{r,s\le0}+\sum_{r>0,\, s\le0} +
\sum_{r\le0,\, s>0}\,.$$ So ergibt sich
$$\eq{\sum_{r,s\le0}{2^{-r-7s/12}\over(2^{-r}+2^{-s})^{29/12}} &\le
2^{-29/12}\sum_{r,s\le0}2^{5r/24+15s/24}\le4\,,\cr \sum_{r>0,\,
s\le0}{2^{-r-7s/12}\over(2^{-r}+2^{-s})^{29/12}}&\le \sum_{r>0,\,
s\le0} 2^{-r}\, 2^{11s/6}\le1.4\,,\cr \sum_{r\le0,\,
s>0}{2^{-r-7s/12}\over(2^{-r}+2^{-s})^{29/12}}&\le \sum_{r\le0,\,
s>0}2^{17r/12}\, 2^{-7s/12}\le3.3\,\cr} $$ und
$$\eq{\sum_{r,s\le0}{2^{-r/2-5s/12}\over(2^{-r}+2^{-s})^{19/12}} &\le
2^{-19/12}\sum_{r,s\le0}2^{7r/24+3s/8}\le8\,,\cr \sum_{r>0,\,
s\le0}{2^{-r/2-5s/12}\over(2^{-r}+2^{-s})^{19/12}}&\le \sum_{r>0,\,
s\le0} 2^{-r/2}\, 2^{7s/6}\le4.4\,,\cr \sum_{r\le0,\,
s>0}{2^{-r/2-5s/12}\over(2^{-r}+2^{-s})^{19/12}}&\le \sum_{r\le0,\,
s>0}2^{13r/12}\, 2^{-5s/12}\le5.7\,.\cr} $$ wobei jeweils in der
ersten Zeile wieder $\sqrt{\al\be}\le(\al+\be)/2$ verwendet wurde.
Insgesamt erhalten wir so
$$ \abs{\sum_{z<g_*({\b m})\le x/y^2}\I {e(g(\b
m)\sqrt{t})\over g^{4}(\b m)}}\le\L^{1/2}\(
1392\,a^{1/6}\,t^{1/12}\,z^{-5/3} +
9123\,a^{1/3}\,t^{-1/12}\,z^{-4/3}\)\,.\eqno(4.19)$$ Verwendet man
(4.18) und (4.19) in (4.8), so ergibt sich
$$ \eq{\abs{P\I(x)}\,\L^{-1/2}&\le\,1213\,a^{1/6}\,(t^*)^{7/12}\,z^{1/3}+4332\,a^{1/3}\,
(t^*)^{5/12}\,z^{2/3}\cr & +
180\,a^{1/6}y^{-2}\,(t^*)^{19/12}\,z^{-5/3} +
1178\,a^{1/3}\,(t^*)^{17/12}y^{-2}\,z^{-4/3}\cr & +
8.4\,x^{1/2}\,y\,\L^{-1/2} + \hbox{ Restterme aus }(4.8)\cr } $$
Durch Ausbalanzieren des ersten gegen den dritten Term auf der
rechten Seite bestimmen wir optimal $z=0.3852\,\sqrt{t^*}/y$ und
erhalten $$ \eq{\abs{P\I(x)}\,\L^{-1/2}&\le\,
1766\,a^{1/6}\,(t^*)^{3/4}\,y^{-1/3} +
6497\,a^{1/3}\,(t^*)^{3/4}\,y^{-1/6}\cr & +
8.4\,x^{1/2}\,y\,\L^{-1/2} + \hbox{ Restterme aus }(4.8)\,.\cr }
$$ Wegen $t^*=x+2y$ und der anf\"{a}nglichen Einschr\"{a}nkung $y\le x/3$
ist $t^*\le5x/3$ und daher $$ \eq{\abs{P\I(x)}\,\L^{-1/2}&\le\,
2591\,a^{1/6}\,x^{3/4}\,y^{-1/3} + 9531\,a^{1/3}\,x^{3/4}\,y^{-2/3}
\cr & + 8.4\,x^{1/2}\,y\,\L^{-1/2} + \hbox{ Restterme aus
}(4.8)\,.\cr }
$$ Balanzieren des ersten gegen den letzten Term rechts ergibt
nun $y=73.6\,a^{1/8}\,x^{3/16}\,\L^{3/8} $ und daher endg\"{u}ltig
$$ \eq{\abs{P\I(x)} &
\le\,1237\,a^{1/8}\,x^{11/16}\,\L^{3/8}+543\,a^{1/4}\,x^{5/8}\,\L^{1/4}\cr
&+\ 268\,{\sqrt{x}\over a}\,\L + 159\,{\sqrt{x}\over
g_0^3}+{2.2\over g_0^3}\,\L+{54.5\over g_0^3}+{50\over
a^{3/2}}\,.\cr}\eqno(4.20)$$ Damit ist der Beweis unseres Satzes
vollendet, soweit es den ersten Fall betrifft, dass der in (4.15)
gew\"{a}hlte Wert f\"{u}r $H$ die Bedingungen (4.12)
erf\"{u}llt.\ssk Wir behandeln nun den Fall jener $N=N_r, M=M_s$,
f\"{u}r die der in (4.15) gegebene Ausdruck gr\"{o}{\ss}er als $\h
M$ ist. Dann ist also $$ M<0.2^{3/4}\,(f(N,M))^{5/4} t^{-1/8}\,.
\eqno(4.21) $$ Wir w\"{a}hlen nun $H:=[\h M]$ (wobei ja wegen (4.7)
$M\ge20$ vorausgesetzt werden kann) und verwenden analog fr\"{u}her
das Superskript $\II$. Aus der Formel vor (4.15) folgt nun $$
{\abs{E_{N,M}(U,W)}\over\sqrt{\log(15N)}} \le 4.7\,M^{1/2}\,N +
14.2\,M^{3/2}N\,f^{-5/4}\,t^{1/8}+ 41.7 \, (MN)^{1/2}
f^{5/4}\,t^{-1/8} $$ anstelle von (4.16), und weiter nach (4.11)
$$ \eq{{\abs{S^{[j]}(N,M)}\over\sqrt{\log(15N)}}\le &\
18.8\,M^{1/2}\,N\,f^{-j} + 56.8\,M^{3/2}N\,f^{-5/4-j}\,t^{1/8}\cr &
+ 166.8\, (MN)^{1/2} f^{5/4-j}\,t^{-1/8}\,.\cr}  $$ Mit (4.21)
ergibt sich $$ \eq{{\abs{S^{[j]}(N,M)}\over\sqrt{\log(15N)}}\le&\
19.6\,N\,f^{5/8-j}\,t^{-1/16}+91.3\,N^{1/2}\,f^{15/8-j}\,t^{-3/16}\cr
\le & {19.6\over
a}\,f^{21/8-j}\,t^{-1/16}+{91.3\over\sqrt{a}}\,f^{23/8-j}\,t^{-3/16}
\cr} \eqno(4.22) $$ wegen $aN\le f^2(N,M)$. F\"{u}r die Teile der
ersten Summe in (4.8) folgt $f(N,M)\le\wz z$ aus $g_*(\b m)\le z$,
und damit $$  \abs{S^{[2]}(N,M)} \L^{-1/2}\le {24.4\over
a}\,z^{5/8}\,t^{-1/16}+{123.7\over\sqrt{a}}\,z^{7/8}\,t^{-3/16}\,.
$$ Wegen $\max(aN,M^2/a^2)\le z^2$, $N\ge10$, $M\ge20$, ist die
Anzahl der betreffenden (dyadisch gew\"{a}hlten) $(N,M)$ durch
$\L^2$ beschr\"{a}nkt. Daher folgt insgesamt
$$ \abs{{\sqrt{t}\over\pi}\sum_{0<g_*({\b m})\le z}\II {e(g(\b m)\sqrt{t})\over g^{2}(\b
m)}}\le\L^{5/2}\({15.6\over
a}\,t^{7/16}\,z^{5/8}+{79.2\over\sqrt{a}}\,t^{5/16}\,z^{7/8}\)\,.
$$ Wegen unserer Wahl $z\le0.3852\sqrt{5x/3}/y$, $y=73.6\,a^{1/8}\,x^{3/16}\,\L^{3/8}
\le x/3$, ergibt dies $$
\eq{\max_{|t-x|\le2y}\abs{{\sqrt{t}\over\pi}\sum_{0<g_*({\b m})\le
z}\II {e(g(\b m)\sqrt{t})\over g^{2}(\b m)}}\le\ &
\,a^{-69/64}\,x^{81/128}\,\L^{145/64}\cr &
+1.2\,a^{-39/64}\,x^{75/128}\,\L^{139/64}\,.\cr }  \eqno(4.23)$$ Wir
verwenden nun (4.22) f\"{u}r $j=4$, um die zweite Summe in (4.8)
abzusch\"{a}tzen. Wegen $f(N,M)\ge g_*(\sqrt{N},0,M)>z$ folgt $$
\abs{S^{[4]}(N,M)} \L^{-1/2}\le {19.6\over
a}\,z^{-11/8}\,t^{-1/16}+{91.3\over\sqrt{a}}\,z^{-9/8}\,t^{-3/16}\,.
$$ Wegen $\max(aN,M^2/a^2)=g_*^2(\sqrt{N},0,M)\le x^2/y^4$ ist die
Zahl der betreffenden dyadi\-schen $(N,M)$ durch $4\L^2$
beschr\"{a}nkt. Daher ergibt sich $$
\abs{{4t^{3/2}\over\pi^3\,y^2}\sum_{z<g_*({\b m})\le x/y^2}\II
{e(g(\b m)\sqrt{t})\over g^{4}(\b m)}}\le\ {\L^{5/2}\over y^2}
\({20.4\over
a}\,z^{-11/8}\,t^{23/16}+{94.4\over\sqrt{a}}\,z^{-9/8}\,t^{21/16}\)\,.$$
Analog fr\"{u}her ben\"{u}tzen wir
$y=73.6\,a^{1/8}\,x^{3/16}\,\L^{3/8} \le x/3$,
$z\ge0.3852\sqrt{x}/y$, und erhalten damit $$
\eq{\max_{|t-x|\le2y}\abs{{4t^{3/2}\over\pi^3\,y^2}\sum_{z<g_*({\b
m})\le x/y^2}\II {e(g(\b m)\sqrt{t})\over g^{4}(\b m)}}\le&\
10.8\,a^{-69/64}\,x^{81/128}\,\L^{145/64}\cr
&+12.6\,a^{-39/64}\,x^{75/128}\,\L^{139/64}\,.\cr} $$ Einsetzen in
(4.8) zusammen mit (4.23) ergibt $$ \abs{P\II(x)}\le
12\(a^{-69/64}\,x^{81/128}\,\L^{145/64} +
a^{-39/64}\,x^{75/128}\,\L^{139/64}\)\,. \eqno(4.24) $$ Wir
behandeln nun zuletzt den Fall (symbolisiert durch $\III$), dass der
Ausdruck f\"{u}r $H$ in (4.15) im Widerspruch zu (4.12) kleiner als
10 ist. Dann ist also wegen $M^2/a^2\le f^2$ $$
f(N,M)<100^{3/4}\,a^{1/4}\,t^{1/8}=: Z^*\,.\eqno(4.25)$$ Wir
w\"{a}hlen nun einfach $H=1$, dann folgt aus der Formel vor (4.15)
$$ \eq{{|E|\over\sqrt{\log(15N)}} &\le
3.1\,{MN}+16.8\,{M^{5/4}\,N\,t^{1/8} \over(f(N,M))^{5/4}}+
34.1\,{M^{3/4}N^{1/2}\, (f(N,M))^{5/4}\over t^{1/8}}\cr & \le
1.2\,f^3 + 5.7\,a^{-1/4}\,f^2\,t^{1/8}+
18.8\,a^{-1/4}\,f^3\,t^{-1/8}\,, \cr }  $$ unter mehrfacher
Verwendung der Mittelungleichung in der Form
$$ M^{2p}\,N^q \le
{p^p\,q^q\over(p+q)^{(p+q)}}\,a^{q-2p}\,(f(N,M))^{2(p+q)}\qquad
(p,q>0)\,.$$ Mit (4.11) und (4.25) erhalten wir
$$ \abs{S^{[2]}(N,M)}\,\L^{-1/2}\le 151\,a^{1/4}\,t^{1/8} +
22.8\,a^{-1/4}\,t^{1/8} + 2400\,.  $$ Wieder ist die Zahl der
betreffenden $(N,M)$ durch $\L^2$ beschr\"{a}nkt, und es folgt
$$ \abs{{\sqrt{t}\over\pi}\sum_{0<g_*({\b m})\le Z^*}\III {e(g(\b m)\sqrt{t})\over g^{2}(\b
m)}}\le{\L^{5/2}\over\pi}\,\(302\,a^{1/4}\,t^{5/8} +
45.6\,a^{-1/4}\,t^{5/8} + 4800\,t^{1/2}\)\,,  $$ und weiter wegen
$y\le x/3$
$$ \max_{|t-x|\le2y}\abs{{\sqrt{t}\over\pi}\sum_{0<g_*({\b m})\le Z^*}\III {e(g(\b m)\sqrt{t})\over g^{2}(\b
m)}}\le{\L^{5/2}}\(134\,a^{1/4}x^{5/8} + 20\,a^{-1/4}x^{5/8} +
2000\,x^{1/2}\). \eqno(4.26) $$ Nach (4.25) ist $Z^*<z$ solange
$|\log a|$ klein gegen $\log x$ ist, der Anteil der zweiten Summe in
(4.8) im Fall $\III$ dann also leer. Um dieses Problem allgemein zu
umgehen, vereinbaren wir, dass f\"{u}r jene $\b m\in\Z^3$, die einem
$(N,M)$-Paar M.d.E.~$\III$ entsprechen, $I(\b m,t)$ stets mittels
(4.5) behandelt wird. Es bleibt daher nur mehr der Beitrag des
zweiten Terms in (4.5) abzusch\"{a}tzen: Wieder nach Lemma 2 folgt
$$ {1\over2\pi^2}\sum_{0<g_*(\b m)\le Z^*} g^{-3}(\b m)\le {21\over\pi^2
 g_0^3}\,\log_+(1.4\,Z^*/g_0)\le {2.2\over g_0^3}\,\L\,. $$ Damit haben
 wir insgesamt gezeigt:
$$ \abs{P\III(x)}\le {\L^{5/2}}\,\(134\,a^{1/4}\,x^{5/8} + 20\,a^{-1/4}\,x^{5/8}
+ 2000\,x^{1/2}\)+ {2.2\over g_0^3}\,\L\,. $$ Zusammen mit (4.20)
und (4.24) ist damit der Beweis unseres Satzes abgeschlossen. \qed

\vbox{\vskip 2.5true cm}

\klein \parindent=0pt

\cen{\bf Literatur}  \bsk \def\smc{}

[1] Bentkus V, G\"otze F (1997) On the lattice point problem for
ellipsoids. Acta Arith {\bf80}: 101--125 \ssk

[2] Chamizo F (1998) Lattice points in bodies of revolution. Acta
Arith {\bf85}: 265-277 \ssk

[3] Chamizo F, Iwaniec H (1995) On the sphere problem. Rev Mat
Iberoamericana {\bf 11}: 417-429 \ssk

[4] G\"otze F (2004) Lattice point problems and values of quadratic
forms. Invent Math {\bf157}: 195-226  \ssk

[5] Heath-Brown DR (1999) Lattice points in the sphere. In: Gy\"ory
et al. (eds.) Number theory in progress, vol.~2, 883-892. Berlin: de
Gruyter \ssk

[6] Hlawka E (1954) \"Uber Integrale auf konvexen K\"orpern I.
Monatsh Math {\bf 54}: 1-36, II, ibid.~{\bf 54}: 81-99 \ssk

[7] Ivi\'c A, Kr\"atzel E, K\"uhleitner M, Nowak WG (2006) Lattice
points in large regions and related arithmetic functions: Recent
developments in a very classic topic. In: Schwarz W (ed.)
Proceedings Conf Elementary and Analytic Number Theory ELAZ'04,
Mainz, im Druck \ssk

[8] Kr\"atzel E (1988) Lattice points. Berlin: VEB Deutscher Verlag
der Wissenschaften \ssk

[9] Kr\"atzel E (2000) Analytische Funktionen in der Zahlentheorie.
Wiesbaden: Teubner.  \ssk

[10] Kr\"atzel E (2004) Lattice points in convex planar domains.
Monatsh Math {\bf143}: 145-162  \ssk

[11] Kr\"atzel E, Nowak WG (2005) Effektive Absch\"{a}tzungen
f\"{u}r den Gitterrest gewisser ebener und dreidimensionaler
Bereiche. Monatsh Math {\bf146}: 21-35 \ssk

[12] Nowak WG (2004) Lattice points in a circle: An improved
mean-square asymptotics. Acta Arith {\bf113}: 259-272 \ssk

[13] Soft Warehouse, Inc. (1995) {\it Derive,} Version 3.11,
Honolulu (Hawaii)  \ssk

[14] Van der Corput JG (1923) Zahlentheoretische Absch\"atzungen mit
Anwendungen auf Gitterpunkt\-probleme. Math.Z. {\bf 17}: 250--259
\ssk

[15] Vinogradov IM (1955) Improvement of asymptotic formulas for the
number of lattice points in a region of three dimensions (Russian).
Izv Akad Nauk SSSR Ser Mat {\bf19}: 3-10 \ssk

[16] Vinogradov IM (1963) On the number of integer points in a
sphere (Russian). Izv Akad Nauk SSSR Ser Mat {\bf27}: 957-968
 \ssk

\vbox{\vskip 1.5true cm}

\parindent=1.5true cm

Adressen der Verfasser: \bsk

\vbox{Ekkehard Kr\"atzel \ssk

Fakult\"{a}t f\"ur Mathematik

Universit\"at Wien

Nordbergstra\ss e 15

1090 Wien, \"Osterreich \ssk

http://www.univie.ac.at/\~{}baxa/kraetzel.html

\bsk

Werner Georg Nowak \ssk

Institut f\"ur Mathematik

Department f\"ur Integrative Biologie

Universit\"at f\"ur Bodenkultur Wien

Gregor Mendel-Stra{\ss}e 33

1180 Wien, \"Osterreich \ssk

E-mail: {\tt \ nowak@mail.boku.ac.at} \ssk

http://www.boku.ac.at/math/nth.html}

\bye